\documentclass[reqno, 12pt]{amsart}
\pdfoutput=1
\makeatletter
\let\origsection=\section \def\section{\@ifstar{\origsection*}{\mysection}}
\def\mysection{\@startsection{section}{1}\z@{.7\linespacing\@plus\linespacing}{.5\linespacing}{\normalfont\scshape\centering\S}}
\makeatother

\usepackage{amsmath,amssymb,amsthm}
\usepackage{mathrsfs}
\usepackage{mathabx}\changenotsign
\usepackage{dsfont}
\usepackage{bbm}

\usepackage{xcolor}
\usepackage[backref=section]{hyperref}
\usepackage[ocgcolorlinks]{ocgx2}
\hypersetup{
	colorlinks=true,
	linkcolor={red!60!black},
	citecolor={green!60!black},
	urlcolor={blue!60!black},
}

\usepackage[open,openlevel=2,atend]{bookmark}

\usepackage[abbrev,msc-links,backrefs]{amsrefs}
\usepackage{doi}

\renewcommand{\PrintDOI}[1]{\doi{#1}}

\usepackage[T1]{fontenc}
\usepackage{lmodern}
\usepackage[babel]{microtype}
\usepackage[english]{babel}

\usepackage{geometry}
\numberwithin{equation}{section}
\numberwithin{figure}{section}

\usepackage{enumitem}

\let\polishlcross=\l
\def\l{\ifmmode\ell\else\polishlcross\fi}

\let\setminus=\smallsetminus

\makeatletter
\def\moverlay{\mathpalette\mov@rlay}
\def\mov@rlay#1#2{\leavevmode\vtop{   \baselineskip\z@skip \lineskiplimit-\maxdimen
		\ialign{\hfil$\m@th#1##$\hfil\cr#2\crcr}}}
\newcommand{\charfusion}[3][\mathord]{
	#1{\ifx#1\mathop\vphantom{#2}\fi
		\mathpalette\mov@rlay{#2\cr#3}
	}
	\ifx#1\mathop\expandafter\displaylimits\fi}
\makeatother

\newcommand{\dcup}{\charfusion[\mathbin]{\cup}{\cdot}}

\DeclareFontFamily{U}  {MnSymbolC}{}
\DeclareSymbolFont{MnSyC}         {U}  {MnSymbolC}{m}{n}
\DeclareFontShape{U}{MnSymbolC}{m}{n}{
	<-6>  MnSymbolC5
	<6-7>  MnSymbolC6
	<7-8>  MnSymbolC7
	<8-9>  MnSymbolC8
	<9-10> MnSymbolC9
	<10-12> MnSymbolC10
	<12->   MnSymbolC12}{}
\DeclareMathSymbol{\powerset}{\mathord}{MnSyC}{180}

\usepackage{tikz}
\usetikzlibrary{calc,decorations.pathmorphing}
\usetikzlibrary{arrows,decorations.pathreplacing}
\pgfdeclarelayer{background}
\pgfdeclarelayer{foreground}
\pgfdeclarelayer{front}
\pgfsetlayers{background,main,foreground,front}

\usepackage{multicol}
\usepackage{subcaption}
\newcommand{\pedge}[9]{
	
	\ifx\relax#6\relax
	\def\qoffs{0pt}
	\else
	\def\qoffs{#6}
	\fi
	
	\def\phedge{
		($#1+#5!\qoffs!-90:#2-#5$) -- 
		($#2+#1!\qoffs!-90:#3-#1$) -- 
		($#3+#2!\qoffs!-90:#4-#2$) -- 
		($#4+#3!\qoffs!-90:#5-#3$) -- 
		($#5+#4!\qoffs!-90:#1-#4$) -- cycle}

	\coordinate (12) at ($#1!\qoffs!90:#2$);
	\coordinate (15) at ($#1!\qoffs!-90:#5$);
	\coordinate (23) at ($#2!\qoffs!90:#3$);
	\coordinate (21) at ($#2!\qoffs!-90:#1$);
	\coordinate (34) at ($#3!\qoffs!90:#4$);
	\coordinate (32) at ($#3!\qoffs!-90:#2$);
	\coordinate (45) at ($#4!\qoffs!90:#5$);
	\coordinate (43) at ($#4!\qoffs!-90:#3$);
	\coordinate (51) at ($#5!\qoffs!90:#1$);
	\coordinate (54) at ($#5!\qoffs!-90:#4$);

	\def\nphedge{
		(15) let \p1=($(15)-#1$), \p2=($(12)-#1$) in 
		arc[start angle={atan2(\y1,\x1)}, delta angle={atan2(\y2,\x2)-atan2(\y1,\x1)-360*(atan2(\y2,\x2)-atan2(\y1,\x1)>0)}, x radius=\qoffs, y radius=\qoffs] --
		(21) let \p1=($(21)-#2$), \p2=($(23)-#2$) in 
		arc[start angle={atan2(\y1,\x1)}, delta angle={atan2(\y2,\x2)-atan2(\y1,\x1)-360*(atan2(\y2,\x2)-atan2(\y1,\x1)>0)}, x radius=\qoffs, y radius=\qoffs] --
		(32) let \p1=($(32)-#3$), \p2=($(34)-#3$) in 
		arc[start angle={atan2(\y1,\x1)}, delta angle={atan2(\y2,\x2)-atan2(\y1,\x1)-360*(atan2(\y2,\x2)-atan2(\y1,\x1)>0)}, x radius=\qoffs, y radius=\qoffs] --
		(43) let \p1=($(43)-#4$), \p2=($(45)-#4$) in 
		arc[start angle={atan2(\y1,\x1)}, delta angle={atan2(\y2,\x2)-atan2(\y1,\x1)-360*(atan2(\y2,\x2)-atan2(\y1,\x1)>0)}, x radius=\qoffs, y radius=\qoffs] --
		(54) let \p1=($(54)-#5$), \p2=($(51)-#5$) in 
		arc[start angle={atan2(\y1,\x1)}, delta angle={atan2(\y2,\x2)-atan2(\y1,\x1)-360*(atan2(\y2,\x2)-atan2(\y1,\x1)>0)}, x radius=\qoffs, y radius=\qoffs] --
		cycle}

	\ifx\relax#7\relax
	\def\plwidth{1pt}
	\else
	\def\plwidth{#7}
	\fi
	
	\ifx\relax#9\relax
	\fill \nphedge;
	\else
	\fill[#9]\nphedge;
	\fi
	
	\ifx\relax#8\relax
	\draw[line width=\plwidth,rounded corners=\qoffs]\nphedge;
	\else
	\draw[line width=\plwidth,#8]\nphedge;
	\fi
}

\newcommand{\qedge}[7]{
	
	\ifx\relax#4\relax
	\def\qoffs{0pt}
	\else
	\def\qoffs{#4}
	\fi
	
	\def\qhedge{
		($#1+#3!\qoffs!-90:#2-#3$) --
		($#2+#1!\qoffs!-90:#3-#1$) --
		($#3+#2!\qoffs!-90:#1-#2$) -- cycle}

	\coordinate (12) at ($#1!\qoffs!90:#2$);
	\coordinate (13) at ($#1!\qoffs!-90:#3$);
	\coordinate (23) at ($#2!\qoffs!90:#3$);
	\coordinate (21) at ($#2!\qoffs!-90:#1$);
	\coordinate (31) at ($#3!\qoffs!90:#1$);
	\coordinate (32) at ($#3!\qoffs!-90:#2$);
	
	\def\nqhedge{
		(13) let \p1=($(13)-#1$), \p2=($(12)-#1$) in
		arc[start angle={atan2(\y1,\x1)}, delta angle={atan2(\y2,\x2)-atan2(\y1,\x1)-360*(atan2(\y2,\x2)-atan2(\y1,\x1)>0)}, x radius=\qoffs, y radius=\qoffs] --
		(21) let \p1=($(21)-#2$), \p2=($(23)-#2$) in
		arc[start angle={atan2(\y1,\x1)}, delta angle={atan2(\y2,\x2)-atan2(\y1,\x1)-360*(atan2(\y2,\x2)-atan2(\y1,\x1)>0)}, x radius=\qoffs, y radius=\qoffs] --
		(32) let \p1=($(32)-#3$), \p2=($(31)-#3$) in
		arc[start angle={atan2(\y1,\x1)}, delta angle={atan2(\y2,\x2)-atan2(\y1,\x1)-360*(atan2(\y2,\x2)-atan2(\y1,\x1)>0)}, x radius=\qoffs, y radius=\qoffs] --
		cycle}
	
	\ifx\relax#5\relax
	\def\qlwidth{1pt}
	\else
	\def\qlwidth{#5}
	\fi
	
	\ifx\relax#7\relax
	\fill \nqhedge;
	\else
	\fill[#7]\nqhedge;
	\fi
	
	\ifx\relax#6\relax
	\draw[line width=\qlwidth,rounded corners=\qoffs]\nqhedge;
	\else
	\draw[line width=\qlwidth,#6]\nqhedge;
	\fi
}

\newcommand{\redge}[8]{
	
	\ifx\relax#5\relax
	\def\qoffs{0pt}
	\else
	\def\qoffs{#5}
	\fi
	
	\def\rhedge{
		($#1+#4!\qoffs!-90:#2-#4$) -- 
		($#2+#1!\qoffs!-90:#3-#1$) -- 
		($#3+#2!\qoffs!-90:#4-#2$) -- 
		($#4+#3!\qoffs!-90:#1-#3$) -- cycle}

	\coordinate (12) at ($#1!\qoffs!90:#2$);
	\coordinate (14) at ($#1!\qoffs!-90:#4$);
	\coordinate (23) at ($#2!\qoffs!90:#3$);
	\coordinate (21) at ($#2!\qoffs!-90:#1$);
	\coordinate (34) at ($#3!\qoffs!90:#4$);
	\coordinate (32) at ($#3!\qoffs!-90:#2$);
	\coordinate (41) at ($#4!\qoffs!90:#1$);
	\coordinate (43) at ($#4!\qoffs!-90:#3$);
	
	\def\nrhedge{
		(14) let \p1=($(14)-#1$), \p2=($(12)-#1$) in 
		arc[start angle={atan2(\y1,\x1)}, delta angle={atan2(\y2,\x2)-atan2(\y1,\x1)-360*(atan2(\y2,\x2)-atan2(\y1,\x1)>0)}, x radius=\qoffs, y radius=\qoffs] --
		(21) let \p1=($(21)-#2$), \p2=($(23)-#2$) in 
		arc[start angle={atan2(\y1,\x1)}, delta angle={atan2(\y2,\x2)-atan2(\y1,\x1)-360*(atan2(\y2,\x2)-atan2(\y1,\x1)>0)}, x radius=\qoffs, y radius=\qoffs] --
		(32) let \p1=($(32)-#3$), \p2=($(34)-#3$) in 
		arc[start angle={atan2(\y1,\x1)}, delta angle={atan2(\y2,\x2)-atan2(\y1,\x1)-360*(atan2(\y2,\x2)-atan2(\y1,\x1)>0)}, x radius=\qoffs, y radius=\qoffs] --
		(43) let \p1=($(43)-#4$), \p2=($(41)-#4$) in 
		arc[start angle={atan2(\y1,\x1)}, delta angle={atan2(\y2,\x2)-atan2(\y1,\x1)-360*(atan2(\y2,\x2)-atan2(\y1,\x1)>0)}, x radius=\qoffs, y radius=\qoffs] --
		cycle}
	
	\ifx\relax#6\relax
	\def\rlwidth{1pt}
	\else
	\def\rlwidth{#6}
	\fi
	
	\ifx\relax#8\relax
	\fill \nrhedge;
	\else
	\fill[#8]\nrhedge;
	\fi
	
	\ifx\relax#7\relax
	\draw[line width=\rlwidth,rounded corners=\qoffs]\nrhedge;
	\else
	\draw[line width=\rlwidth,#7]\nrhedge;
	\fi
}

\let\epsilon=\varepsilon

\let\rho=\varrho
\let\theta=\vartheta

\newtheoremstyle{note}  {4pt}  {4pt}  {\sl}  {}  {\bfseries}  {.}  {.5em}          {}
\newtheoremstyle{introthms}  {3pt}  {3pt}  {\itshape}  {}  {\bfseries}  {.}  {.5em}          {\thmnote{#3}}
\newtheoremstyle{remark}  {2pt}  {2pt}  {\rm}  {}  {\bfseries}  {.}  {.3em}          {}

\theoremstyle{plain}
\newtheorem{theorem}{Theorem}[section]

\theoremstyle{note}

\theoremstyle{remark}

\usepackage{lineno}
\newcommand*\patchAmsMathEnvironmentForLineno[1]{
	\expandafter\let\csname old#1\expandafter\endcsname\csname #1\endcsname
	\expandafter\let\csname oldend#1\expandafter\endcsname\csname end#1\endcsname
	\renewenvironment{#1}
	{\linenomath\csname old#1\endcsname}
	{\csname oldend#1\endcsname\endlinenomath}}
\newcommand*\patchBothAmsMathEnvironmentsForLineno[1]{
	\patchAmsMathEnvironmentForLineno{#1}
	\patchAmsMathEnvironmentForLineno{#1*}}
\AtBeginDocument{
	\patchBothAmsMathEnvironmentsForLineno{equation}
	\patchBothAmsMathEnvironmentsForLineno{align}
	\patchBothAmsMathEnvironmentsForLineno{flalign}
	\patchBothAmsMathEnvironmentsForLineno{alignat}
	\patchBothAmsMathEnvironmentsForLineno{gather}
	\patchBothAmsMathEnvironmentsForLineno{multline}
}

\def\ex{\text{\rm ex}}

\usepackage{scalerel}

\makeatletter
\newcommand{\overrighharpoonup}[1]{\ThisStyle{%
		\vbox {\m@th\ialign{##\crcr
				\rightharpoonupfill \crcr
				\noalign{\kern-\p@\nointerlineskip}
				$\hfil\SavedStyle#1\hfil$\crcr}}}}

\def\rightharpoonupfill{%
	$\SavedStyle\m@th\mkern+0.8mu\cleaders\hbox{$\shortbar\mkern-4mu$}\hfill\rightharpoonuptip\mkern+0.8mu$}

\def\rightharpoonuptip{%
	\raisebox{\z@}[2pt][1pt]{\scalebox{0.55}{$\SavedStyle\rightharpoonup$}}}

\def\shortbar{%
	\smash{\scalebox{0.55}{$\SavedStyle\relbar$}}}
\makeatother

\makeatletter
\newcommand{\overlefharpoonup}[1]{\ThisStyle{%
		\vbox {\m@th\ialign{##\crcr
				\leftharpoonupfill \crcr
				\noalign{\kern-\p@\nointerlineskip}
				$\hfil\SavedStyle#1\hfil$\crcr}}}}

\def\leftharpoonupfill{%
	$\SavedStyle\m@th\mkern+0.8mu\cleaders\hbox{$\shortbar\mkern-4mu$}\hfill\leftharpoonuptip\mkern+0.8mu$}

\def\leftharpoonuptip{%
	\raisebox{\z@}[2pt][1pt]{\scalebox{0.55}{$\SavedStyle\leftharpoonup$}}}

\makeatletter
\newsavebox\myboxA
\newsavebox\myboxB
\newlength\mylenA

\newcommand*\xoverline[2][0.75]{%
	\sbox{\myboxA}{$\m@th#2$}%
	\setbox\myboxB\null
	\ht\myboxB=\ht\myboxA%
	\dp\myboxB=\dp\myboxA%
	\wd\myboxB=#1\wd\myboxA
	\sbox\myboxB{$\m@th\overline{\copy\myboxB}$}
	\setlength\mylenA{\the\wd\myboxA}
	\addtolength\mylenA{-\the\wd\myboxB}%
	\ifdim\wd\myboxB<\wd\myboxA%
	\rlap{\hskip 0.5\mylenA\usebox\myboxB}{\usebox\myboxA}%
	\else
	\hskip -0.5\mylenA\rlap{\usebox\myboxA}{\hskip 0.5\mylenA\usebox\myboxB}%
	\fi}
\makeatother

\begin{document}
	
	\title[Hypergraphs accumulate infinitely often]
	{Hypergraphs accumulate infinitely often}
	
	\author[D.~Conlon]{David Conlon}
	\address{Department of Mathematics, California Institute of Technology, Pasadena, USA}
	\email{dconlon@caltech.edu}
	
	\author[B. Sch\"ulke]{Bjarne Sch\"ulke}
    \address{Extremal Combinatorics and Probability Group, Institute for Basic Science, Daejeon, South Korea}
    \email{schuelke@ibs.re.kr}

	
	\begin{abstract}
        We show that the set $\Pi^{(k)}$ of Tur\'an densities of $k$-uniform hypergraphs has infinitely many accumulation points in $[0,1)$ for every $k \geq 3$. This extends an earlier result of ours showing that $\Pi^{(k)}$ has at least one such accumulation point.
	\end{abstract}
	
	\maketitle
	
	\section{Introduction}
        For~$k\in\mathds{N}$, a \emph{$k$-uniform hypergraph} (or \emph{$k$-graph})~$H=(V,E)$ consists of a vertex set~$V$ and an edge set~$E\subseteq V^{(k)}=\{e\subseteq V:\vert e\vert=k\}$.
        Given~$n\in\mathds{N}$ and a family of $k$-graphs~$\mathcal{F}$, the \emph{extremal number}~$\ex(n,\mathcal{F})$ is the maximum number of edges in a~$k$-graph~$H$ with~$n$ vertices that does not contain a copy of any graph in~$\mathcal{F}$.
        The \emph{Tur\'an density} of~$\mathcal{F}$ is then given by
        $$\pi(\mathcal{F})=\lim_{n\to\infty}\frac{\ex(n,\mathcal{F})}{\binom{n}{k}},$$
        where the limit is known, by a simple monotonicity argument~\cite{KNS:64}, to be well-defined.
        If~$\mathcal{F}=\{F\}$ for some~$k$-graph~$F$, we omit the parentheses, writing~$\pi(\mathcal{F})=\pi(F)$.
        The problem of determining these Tur\'an densities is one of the oldest and most fundamental questions in extremal combinatorics. 
        
        When $k=2$, that is, when $\mathcal{F}$ is a family of graphs,~$\pi(\mathcal{F})$ is essentially completely understood, with the final result, the culmination of work by Tur\'an~\cite{T:41}, Erd\H{o}s and Stone~\cite{ES:46}, and Erd\H{o}s and Simonovits~\cite{ES:66}, saying that~$\pi(\mathcal{F})=\frac{\chi(\mathcal{F})-2}{\chi(\mathcal{F})-1}$, where~$\chi(\mathcal{F})$ is the minimum chromatic number of an element of~$\mathcal{F}$. If we set 
        \begin{align*}
        \Pi^{(k)} & =\{\pi(F):\:F\text{ is a }k\text{-graph}\}\,,\\
        \Pi^{(k)}_{\text{fin}} & =\{\pi(\mathcal{F}):\:\mathcal{F}\text{ is a finite family of }k\text{-graphs}\}\,,\\
        \Pi_{\infty}^{(k)} & =\{\pi(\mathcal{F}):\:\mathcal{F}\text{ is a family of }k\text{-graphs}\}\,,
        \end{align*}
        then the Erd\H{o}s--Stone--Simonovits theorem implies that
        \[\Pi^{(2)}=\Pi^{(2)}_{\text{fin}}=\Pi_{\infty}^{(2)}=\{0, 1/2, 2/3, 3/4, \dots\}.\]
        In particular, each of these sets is well-ordered.

	For~$k \ge 3$, there is no clear analogue of the Erd\H{o}s--Stone--Simonovits theorem for any of the sets~$\Pi^{(k)}$,~$\Pi^{(k)}_{\text{fin}}$, or~$\Pi_{\infty}^{(k)}$.
    For~$\Pi_{\infty}^{(k)}$, this was shown by Frankl and R\"odl~\cite{FR:84}, who proved that~$\Pi^{(k)}_{\infty}$ is not well-ordered.
    That is, the set has downward accumulation points. 
    In particular, this disproved the jumping conjecture of Erd\H{o}s, which suggested that, like for~$\Pi^{(2)}_\infty$, there should be a non-trivial gap or jump between any element of~$[0,1)$ and the next element of~$\Pi^{(k)}_{\infty}$.
    By applying the important result of Pikhurko~\cite{P:14} that~$\Pi_{\infty}^{(k)}$ is the closure of~$\Pi^{(k)}_{\text{fin}}$, it is easily seen that~$\Pi^{(k)}_{\text{fin}}$ is also not well-ordered.

	While it remains an intriguing open problem to show that $\Pi^{(k)}$ is again not well-ordered, a first step showing that $\Pi^{(k)}$ is indeed more complex than $\Pi^{(2)}$ was taken in the recent paper~\cite{CS:25}, where we proved the following result.
	
	\begin{theorem}\label{thm:old}
            For every integer~$k\geq3$, the set~$\Pi^{(k)}$ has an accumulation point in~$[0,1)$.
        \end{theorem}
        
        Here we extend this result, showing that $\Pi^{(k)}$ has infinitely many accumulation points.
        This goes another step further in showing how much more complex~$\Pi^{(k)}$ gets for~$k\geq3$.
        
        \begin{theorem}\label{thm:main}
            For every integer~$k\geq3$, the set~$\Pi^{(k)}$ has infinitely many accumulation points in~$[0,1)$.
        \end{theorem}

        This is a consequence of the following result which states that, in addition, for each of these accumulation points~$\alpha$ there is a family of~$k$-graphs whose Tur\'an density is~$\alpha$.

        \begin{theorem}\label{thm:extended}
            For every integer~$k\geq3$, there are infinitely many~$\alpha\in [0,1)$ such that there are two sequences of~$k$-graphs,~$\{F_i\}_{i\in\mathds{N}}$ and~$\{G_i\}_{i\in\mathds{N}}$, with the following properties:
            \begin{enumerate}
                \item[I.]\label{it:thm:accpoint} $\pi(F_i)\to\alpha$ and~$\pi(F_i)<\alpha$ for all~$i\in\mathds{N}$.
                \item[II.]\label{it:thm:lim} For all~$\varepsilon>0$, there is some~$i\in\mathds{N}$ such that~$\alpha\leq\pi(G_i)\leq\alpha+\varepsilon$.
             \end{enumerate}
        \end{theorem}

	The proof of this result will occupy the remainder of this short paper.
        
    \section{Preliminaries} \label{sec:pre}
        Given an integer~$t$ and a~$k$-graph~$F$, let~$B(F,t)$ be the \emph{$t$-blow-up} of~$F$, the~$k$-graph obtained from~$F$ by replacing every vertex by~$t$ copies of itself.
        The following phenomenon, which we make extensive use of, is well-known (see, for instance, Lemma~2.1 and Theorem~2.2 in~\cite{K:11}, as well as the subsequent discussion).
        
        \begin{theorem}[Supersaturation]\label{thm:supersaturation}
            \begin{enumerate}
                \item\label{it:supsat:original} For every~$k$-graph~$F$ and~$\delta>0$, there are~$\varepsilon>0$ and~$n_0$ such that every~$k$-graph on~$n\geq n_0$ vertices with at least~$(\pi(F)+\delta)\binom{n}{k}$ edges contains at least~$\varepsilon n^{\vert V(F)\vert}$ copies of~$F$.
                \item\label{it:supsat:blowup} For every integer~$t$ and~$k$-graph~$F$, $\pi(B(F,t))=\pi(F)$.
                \item\label{it:supsat:hom} Let~$F$ be a~$k$-graph and $\mathcal{F}$ be the (finite) family of~$k$-graphs~$F'$ whose vertex set is a subset of~$V(F)$ and for which there exists a homomorphism~$\varphi:F\to F'$.
                Then~$\pi(\mathcal{F})=\pi(F)$.
                \item\label{it:supsat:onevtx} For every~$k$-graph~$F$ and~$\delta>0$, there are~$\varepsilon>0$ and~$n_0$ such that, for all~$v\in V(F)$, every~$k$-graph on~$n\geq n_0$ vertices with at least~$(\pi(F)+\delta)\binom{n}{k}$ edges contains the~$k$-graph obtained from~$F$ by replacing~$v$ by~$\varepsilon n$ copies of~$v$.
            \end{enumerate}
        \end{theorem}

        We will also make use of expansions of hypergraphs.
        Setting~$X=\{x_1,\dots,x_s\}$, the~$k$-uniform \emph{expansion} of~$K_s^{(2)}$, the complete~$2$-graph on~$s$ vertices, is the~$k$-graph~$G_s^{(k)}$ with vertex set 
        $$X\dcup\{v_i^e:\:i\in[k-2],\,e\in X^{(2)}\}$$ 
        and edge set
        $$\{e\cup\{v_1^e,\dots,v_{k-2}^e\}:\:e\in X^{(2)}\}.$$
        In other words, the~$k$-uniform expansion of~$K_s^{(2)}$ is obtained from~$K_s^{(2)}$ by adding~$k-2$ new vertices to each edge.
        We will need the following result of Mubayi~\cite{M:06} determining the Tur\'an density of these expansions.

        \begin{theorem}\label{thm:expansion}
            For all integers~$s> k\geq2$,~$\pi(G_{s}^{(k)})=\frac{(s-1)\cdot(s-2)\cdots(s-k)}{(s-1)^k}$.
        \end{theorem}
    
    \section{Proof of Theorem~\ref{thm:extended}}
        The proof makes use of some $k$-graphs that are obtained by gluing ladders and zycles, both of which we now define, in appropriate ways.
        
        For~$k,\ell\in\mathds{N}$, we define the~$k$-uniform \emph{ladder of length~$\ell$} to be the~$k$-graph~$L^{(k)}_{\ell}$ with vertex set $$V(L^{(k)}_{\ell})=\{v_{ij}:\:i\in[\ell],\,j\in[k-1]\}\dcup\{t\}$$ and edge set $$E(L^{(k)}_{\ell})=\{v_{i1}\dots v_{ik-1}v_{i+1j}:\:i\in[\ell-1],\,j\in[k-1]\}\cup\{v_{\ell1}\dots v_{\ell k-1}t\}\,.$$

        For $k,\ell\in\mathds{N}$ with~$\ell\geq 2$, we define the $k$-uniform \emph{zycle of length~$\ell$} to be the~$k$-graph~$Z^{(k)}_{\ell}$ with vertex set
        $$V(Z^{(k)}_{\ell})=\{v_{ij}:i\in\mathds{Z}/\ell\mathds{Z},j\in[k-1]\}$$ 
        and edge set
        $$E(Z^{(k)}_{\ell})=\{v_{i1}\dots v_{ik-1}v_{i+1j}:\:i\in\mathds{Z}/\ell\mathds{Z},\,j\in[k-1]\}\,.$$

        Next we define a~$k$-graph that, roughly speaking, is obtained by gluing a zycle of length~$\ell$ to the last~$(k-1)$-set of a ladder of length~$m$. 
        More formally, for $k,m,\ell\in\mathds{N}$ with $\ell \ge 2$, the~$k$-graph~$LZ^{(k)}(m,\ell)$ is the~$k$-graph with vertex set
        \begin{align*}
            V(LZ^{(k)}(m,\ell)) = \,
            & \{v_{ij}:\:i\in[m],\,j\in[k-1]\}\\
            & \dcup\{w_{ij}:\:i\in[\ell]\setminus\{1\},\,j\in[k-1]\}
        \end{align*}
        and edge set
        \begin{align*}
            E(LZ^{(k)}(m,\ell))= \,
            &\{v_{i1}\dots v_{ik-1} v_{i+1j}:\:i\in[m-1],\,j\in[k-1]\}\\
            &\cup\{v_{m1}\dots v_{m k-1} w_{2j}:\:j\in[k-1]\}\\
            &\cup\{w_{i1}\dots w_{ik-1} w_{i+1j}:\:i\in[\ell-1]\setminus\{1\},\,j\in[k-1]\}\\
            &\cup\{w_{\ell1}\dots w_{\ell k-1} v_{mj}:\:j\in[k-1]\}\,.
        \end{align*}
        For each of~$L^{(k)}_{\ell}$,~$Z^{(k)}_{\ell}$, and~$LZ^{(k)}(m,\ell)$, we sometimes refer to the set~$\{v_{11},\dots,v_{1k-1}\}$ as the \emph{starting set}.
        
        If, in addition,~$s\in\mathds{N}$ with~$s>k$, we consider the~$k$-graph obtained from an~$s$-set~$X$ by adding, for each pair~$e\in X^{(2)}$, a copy of~$L_{\ell_{e}}^{(k)}$ such that these copies only intersect in vertices of~$X$, 
        where the length~$\ell_e$ may depend on the pair~$e$.
        More formally, given a set~$X=\{x_1,\dots,x_s\}$, let~$\{e_1,\dots,e_{\binom{s}{2}}\}=X^{(2)}$ be an enumeration of the pairs of elements in~$X$.
        Furthermore, let~$\ell_1,\dots,\ell_{\binom{s}{2}}\in\mathds{N}$ with~$\ell_i\geq\ell_{i+1}$ for all~$i\in[\binom{s}{2}-1]$. 
        For~$x_ix_j=e\in X^{(2)}$ with~$i<j$, we write~$v_{1k-2}^e=x_i$ and~$v_{1k-1}^e=x_j$.
        For all~$e\in X^{(2)}$ and $j\in[k-3]$, let~$v_{1j}^{e}$ be pairwise distinct vertices which are also distinct from any vertex in~$X$.
        Finally, for all~$e_r\in X^{(2)}$,~$i\in[\ell_r]\setminus\{1\}$, and~$j\in[k-1]$, let~$v_{ij}^{e_r}$ be distinct vertices which are also distinct from any previously chosen vertices.
        Then we define~$GL^{(k)}(s;\ell_1,\dots,\ell_{\binom{s}{2}})$ to be the~$k$-graph with vertex set
        \begin{align}\label{eq:vtxsetGL}
            V(GL^{(k)}(s;\ell_1,\dots,\ell_{\binom{s}{2}}))= \,
            &\{v_{ij}^{e_r}:\:e_r\in X^{(2)},\,i\in[\ell_r],\,j\in[k-1]\}\nonumber\\
            &\dcup\{t^e:e\in X^{(2)}\}
        \end{align}
        and edge set
        \begin{align*}
            E(GL^{(k)}(s;\ell_1,\dots,\ell_{\binom{s}{2}}))= \,
            &\{v_{i1}^{e_r}\dots v_{ik-1}^{e_r} v_{i+1j}^{e_r}:\:e_r\in X^{(2)},\,i\in[\ell_r-1],\,j\in[k-1]\}\\
            &\cup\{v_{\ell_r 1}^{e_r}\dots v_{\ell_r k-1}^{e_r} t^{e_r}:\:e_r\in X^{(2)}\}\,.
        \end{align*}
        If~$\ell_i=\ell$ for all~$i\in[\binom{s}{2}]$, we simply write~$GL^{(k)}(s,\ell)$ for~$GL^{(k)}(s;\ell,\dots,\ell)$.
        We also note that~$GL^{(k)}(s,1)=G_s^{(k)}$, as defined in Section~\ref{sec:pre}.

        Lastly, we need one more type of~$k$-graph.
        Roughly speaking, it is obtained from $GL^{(k)}(s,\ell)$ by closing the ends of the ladders that come out of the set~$X$ into zycles.
        More formally, for $k,s,m,\ell\in\mathds{N}$ with $s > k$ and $\ell \ge 2$, we define the~$k$-graph~$GLZ^{(k)}(s,m,\ell)$ to be the~$k$-graph with vertex set
        \begin{align*}
            V(GLZ^{(k)}(s,m,\ell))= \,
            &\{v_{ij}^e:\:e\in X^{(2)},\,i\in[m],\,j\in[k-1]\}\\
            &\dcup\{w^e_{ij}:\:e\in X^{(2)},\,i\in[\ell]\setminus\{1\},\,j\in[k-1]\}
        \end{align*}
        and edge set
        \begin{align*}
            E(GLZ^{(k)}(s,m,\ell))=\,
            &\{v_{i1}^e\dots v_{ik-1}^e v_{i+1j}^e:\:e\in X^{(2)},\,i\in[m-1],\,j\in[k-1]\}\\
            &\cup\{v_{m1}^e\dots v_{m k-1}^e w_{2j}^e:\:e\in X^{(2)},\,j\in[k-1]\}\\
            &\cup\{w_{i1}^e\dots w_{ik-1}^e w_{i+1j}^e:\:e\in X^{(2)},\,i\in[\ell-1]\setminus\{1\},\,j\in[k-1]\}\\
            &\cup\{w_{\ell1}^e\dots w_{\ell k-1}^e v_{mj}^e:\:e\in X^{(2)},\,j\in[k-1]\}\,.
        \end{align*}
        
        From now on, we suppress the uniformity in the notation if it is clear from context, for instance, writing~$GL(s,\ell)$ instead of~$GL^{(k)}(s,\ell)$. In outline, the proof of Theorem~\ref{thm:extended} will proceed as follows. 
        First, we show that for every integer~$s>k \ge 3$ there is some~$\alpha_s\in[0,1]$ such that~$\lim_{\ell\to\infty}\pi(GL(s,\ell))=\alpha_s$,  but~$\pi(GL(s,\ell))<\alpha_s$ for all~$\ell\in\mathds{N}$.
        Because~$GL(s,\ell)\subseteq GL(s+1,\ell)$ for all~$s>k$ and $\ell$, we have~$\alpha_s\leq\alpha_{s+1}$.
        We will argue that, more strongly, for every~$s>k$ there is some~$s'$ such that~$\alpha_s<\alpha_{s'}$.
        Together, these imply Part I of the theorem.
        We will then show that for every~$s>k$ and~$\varepsilon>0$, there are~$m, \ell\in\mathds{N}$ such that~$\alpha_s\leq\pi(GLZ(s,m,\ell))\leq\alpha_s+\varepsilon$, which will complete the proof.
        
        \subsection{Part I}
        Let~$s>k \ge 3$ be an integer.
        To show that there is some~$\alpha_s\in[0,1]$ such that~$\lim_{\ell\to\infty}\pi(GL(s,\ell))=\alpha_s$, but~$\pi(GL(s,\ell))<\alpha_s$ for all~$\ell\in\mathds{N}$, it is sufficient to show that~$\pi(GL(s,\ell))<\pi(GL(s,\ell+1))$ for all~$\ell\in\mathds{N}$.
        We do this by induction on~$\ell$. 
        
        First, let~$\ell=1$ and, for a given~$n\in\mathds{N}$, let~$H$ be the~$k$-graph that is obtained from a balanced complete~$(s-1)$-partite~$k$-graph with partition~$V(H)=[n]=V_1\dcup\cdots\dcup V_{s-1}$ by adding a balanced complete~$k$-partite~$k$-graph inside each partition class (with partition~$V_i=W_i^1\dcup\cdots\dcup W_i^k$ for each~$i\in[s-1]$).
        If there were a copy of~$GL(s,2)$ in~$H$, there would have to be at least two vertices~$x,x'\in X$ that lie in the same partition class~$V_i$.
        By the constructions of~$GL(s,2)$ and~$H$, the vertices~$v_{21}^{xx'}$ and~$v_{22}^{xx'}$, say, must lie in the same~$W_i^j$. 
        But then in~$H$ there is no edge containing both~$v_{21}^{xx'}$ and~$v_{22}^{xx'}$ (which exists in~$GL(s,2)$), meaning that, in fact,~$H$ has to be~$GL(s,2)$-free.
        This implies that~$\pi(GL(s,2))>k!\binom{s-1}{k}\frac{1}{(s-1)^k}=\pi(GL(s,1))$, where the last inequality comes from  Theorem~\ref{thm:expansion}.
        
        Now assume that~$\ell>1$ and that~$\pi(GL(s,i))<\pi(GL(s,i+1))$ holds for all~$i\in[\ell-1]$.
        We will show that~$\pi(GL(s,\ell))<\pi(GL(s,\ell+1))$.

        By induction, we know that~$\pi(GL(s,\ell))>\pi(GL(s,\ell-1))$.
        Thus, there is some maximum~$r\in[\binom{s}{2}]$ such that, setting~$\ell_i=\ell$ for~$i\in[r-1]$ and~$\ell_i=\ell-1$ for~$i\in[r,\binom{s}{2}]$, we have~$\pi(GL(s;\ell_1,\dots,\ell_{\binom{s}{2}}))<\pi(GL(s,\ell))$.
        Let~$\ell'_i=\ell$ for~$i\in[r]$ and~$\ell'_i=\ell-1$ for~$i\in[r+1,\binom{s}{2}]$. Denote by~$\mathcal{GL}$ the (finite) family of~$k$-graphs~$F$ whose vertex set is a subset of~$V(GL(s;\ell'_1,\dots,\ell'_{\binom{s}{2}}))$ and for which there exists a homomorphism~$\varphi:\:GL(s;\ell'_1,\dots,\ell'_{\binom{s}{2}})\to F$.
        By supersaturation (Theorem~\ref{thm:supersaturation}~\eqref{it:supsat:hom}) (and the choice of~$r$), we know that~$\pi(\mathcal{GL})=\pi(GL(s,\ell))$ and thus it suffices to show that~$\pi(\mathcal{GL})<\pi(GL(s,\ell+1))$.

        Set~$\pi_0=\pi(GL(s;\ell_1,\dots,\ell_{\binom{s}{2}}))$ and note that~$\pi(\mathcal{GL})>\pi_0$.
        Therefore, setting~$\eta=\pi(\mathcal{GL})-\pi_0$, we have~$\eta>0$.
        Furthermore, by supersaturation (Theorem~\ref{thm:supersaturation}~\eqref{it:supsat:onevtx}), there is some $\varepsilon_1>0$ such that, for~$n$ sufficiently large, every~$k$-graph~$H$ on~$n$ vertices with at least~$(\pi_0+\eta/2)\binom{n}{k}$ edges contains a copy of the~$k$-graph~$G$ that is obtained from~$GL(s;\ell_1,\dots,\ell_{\binom{s}{2}})$ by blowing up the vertex~$t^{e_r}$ to a set~$T$ of size~$\varepsilon_1n$.
        Finally, let~$\varepsilon_2\ll\varepsilon_1,\eta$ with~$\varepsilon_2>0$ and let~$n\in\mathds{N}$ be sufficiently large that\footnote{By the monotonicity argument mentioned in the introduction, all of the terms on the left-hand side are non-negative.}
        \begin{align}\label{eq:constants}
            &\frac{\ex(n,GL(s;\ell_1,\dots,\ell_{\binom{s}{2}}))}{\binom{n}{k}}-\pi_0<\varepsilon_2\,,\nonumber\\
            &\frac{\ex(n,\mathcal{GL})}{\binom{n}{k}}-\pi(\mathcal{GL})<\varepsilon_2\,,\,\text{and}\nonumber\\
            &\frac{\ex(n,GL(s,\ell+1))}{\binom{n}{k}}-\pi(GL(s,\ell+1))<\varepsilon_2\,.
        \end{align}

        Now consider an extremal example~$H$ for~$\mathcal{GL}$ on~$n$ vertices.
        By our choice of constants, we know that~$H$ contains a copy of~$G$.
        If any~$(k-1)$-subset of~$T$ is contained in an edge of~$H$, then~$H$ would contain a (possibly) degenerate copy of~$GL(s;\ell'_1,\dots,\ell'_{\binom{s}{2}})$, i.e., a copy of an element of~$\mathcal{GL}$.
        Thus, no~$(k-1)$-subset of~$T$ is contained in an edge of~$H$.
        
        Next we add to~$H$ a complete balanced $k$-partite~$k$-graph on~$T=T_1\dcup\cdots\dcup T_k$ and call the resulting~$k$-graph~$H'$.
        We claim that~$H'$ is~$GL(s,\ell+1)$-free.
        Assume, for the sake of contradiction, that~$H'$ contains a copy of~$GL(s,\ell+1)$ with vertex set as in~\eqref{eq:vtxsetGL}.\footnote{To avoid making the notation messier, we will not give the vertices new names. We do not mean the vertices of this copy of $GL(s,\ell+1)$ to be necessarily the same as some of the vertices of the copy of $G$.}
        Since this copy of~$GL(s,\ell+1)$ is not contained in~$H$, one of its edges must be an edge~$z_1\dots z_k\in E(H')\setminus E(H)$, so we also have~$z_1,\dots, z_k\in T$.
        In fact, since~$H$ is (in particular)~$GL(s,\ell)$-free, there must be~$e\in X^{(2)}$,~$i\in[\ell]$, and~$j\in[k-1]$ such that~$z_1\dots z_k$ is one of the edges~$v^e_{i1}\dots v^e_{ik-1}v^e_{i+1j}$.
        Without loss of generality, assume that~$z_1=v^e_{i1},\dots,z_{k-1}=v^e_{ik-1}$ for some~$e\in X^{(2)}$,~$i\in[\ell]$, and~$j\in[k-1]$ with~$v^e_{i1}\in T_1,\dots,v^e_{ik-1}\in T_{k-1}$.
        Recall that, by the construction of~$H'$ and the discussion in the previous paragraph, any edge of~$H'$ containing a~$(k-1)$-subset of~$T$ must be in~$E(H')\setminus E(H)$ and must therefore contain exactly one vertex from each of~$T_1,\dots,T_k$.
        Thus,~$v^e_{(i+1)1},\dots,v^e_{(i+1)k-1}\in T_k$ and so these~$k-1$ vertices cannot lie together in any edge of~$H'$, contradicting that there is 
        a copy of~$GL(s,\ell+1)$ in~$H'$.
        Hence,~$H'$ is indeed a~$GL(s,\ell+1)$-free~$k$-graph on~$n$ vertices.
        
        By monotonicity, we know that~$H$ has at least~$\pi(\mathcal{GL})\binom{n}{k}$ edges.
        Therefore,~$H'$ has more than 
        $$\pi(\mathcal{GL})\binom{n}{k}+\left(\frac{\varepsilon_1n}{k+1}\right)^k>(\pi(\mathcal{GL})+\varepsilon_2)\binom{n}{k}$$ 
        edges.
        By~\eqref{eq:constants}, this means that~$\pi(GL(s,\ell+1))>\pi(\mathcal{GL})=\pi(GL(s,\ell))$.
        We have therefore proved that~$\lim_{\ell\to\infty}\pi(GL(s,\ell))=\alpha_s$  for some~$\alpha_s\in[0,1]$ with~$\pi(GL(s,\ell))<\alpha_s$ for all~$\ell\in\mathds{N}$.

        Next we argue that for every integer~$s>k$, there is some integer~$s'\gg s$ such that~$\alpha_{s'} > \alpha_s$.
        Note that since~$GL(s',\ell)\supseteq G_{s'}$ for every $\ell \in \mathds{N}$, Theorem~\ref{thm:expansion} implies that  $$\alpha_{s'}\geq\frac{(s'-1)\cdot(s'-2)\cdots(s'-k)}{(s'-1)^k}\,.$$
        On the other hand, observe that~$GL(s,\ell)$ is contained in a blow-up of~$K^{(k)}_{s+\binom{s}{2}\cdot[(k-3)+(k-1)]}$ for every~$\ell\in\mathds{N}$.
        Therefore, by Theorem~\ref{thm:supersaturation}~\eqref{it:supsat:blowup},~$\alpha_s\leq\pi(K^{(k)}_{s+\binom{s}{2}\cdot(2k-4)})$.
        Since $$\frac{(s'-1)\cdot(s'-2)\cdots(s'-k)}{(s'-1)^k}\to 1$$ as~$s'\to\infty$ and~$\pi(K^{(k)}_{s+\binom{s}{2}\cdot(2k-4)})<1$, we indeed have~$\alpha_{s'}>\alpha_s$ for~$s'\gg s$.

        \subsection{Part II}
        Let~$s>k$ be an integer and let~$\varepsilon>0$.
        Choose~$t,n\in\mathds{N}$ such that $$\varepsilon,s^{-1}\gg t^{-1}\gg n^{-1}$$ and, for simplicity, assume that~$t\mid n$.
        Now let~$H$ be a~$k$-graph with vertex set~$[t]$ and~$e(H)\geq(\alpha_s+\varepsilon)\binom{t}{k}$.
        We will show that there is a homomorphism from~$GLZ(s,\binom{t}{k-1}+1,\binom{t}{k-1}!)$ into~$H$.
        Let~$H_*=B(H,n/t)$ be the~$k$-graph obtained from~$H$ by replacing every vertex~$i$ of~$H$ by~$n/t$ copies of itself, the set of which we call~$V_i$.
        For~$v\in V(H_*)$, let~$f(v)$ denote the index of the partition class of~$H_*$ that contains~$v$, i.e., if~$v$ is one of the copies of the vertex~$i\in V(H)$, then~$f(v)=i$.
        Then~$H_*$ is a~$k$-graph on~$n$ vertices with~$e(H_*)\geq(\alpha_s+\varepsilon)\binom{t}{k}\big(\frac{n}{t}\big)^k\geq(\alpha_s+\varepsilon/2)\binom{n}{k}$.
        
        Since~$\pi(GL(s,\ell))<\alpha_s$ for every~$\ell\in\mathds{N}$, we have that~$H_*$ contains a copy of $GL(s,\binom{t}{k-1}+1)$ with vertex set as in~\eqref{eq:vtxsetGL}.
        Fix~$e\in X^{(2)}$.
        Note that for each~$i\in[\binom{t}{k-1}+1]$, the indices~$f(v^e_{ij})$ with~$j\in[k-1]$ are pairwise distinct, since~$v^e_{i1},\dots,v^e_{ik-1}$ are contained in an edge together.
        As~$H_*$ only has~$t$ distinct partition classes, we deduce from the pigeonhole principle that, for some~$i,i'\in[\binom{t}{k-1}+1]$ with~$i'<i$, we have $$\{f(v^e_{i1}),\dots,f(v^e_{ik-1})\}=\{f(v^e_{i'1}),\dots,f(v^e_{i'k-1})\}\,.$$
        Since~$H_*$ is a blow-up of~$H$, this implies that there is a homomorphism of a zycle of length at most~$\binom{t}{k-1}$ into~$H$ that maps the starting set to~$\{f(v^e_{i1}),\dots,f(v^e_{ik-1})\}$.
        As described in~\cite{PS:23}, ``cycling'' through any such zycle the right number of times yields a homomorphism from~$Z_{\binom{t}{k-1}!}$ to~$H$ that maps the starting set to~$\{f(v^e_{i1}),\dots,f(v^e_{ik-1})\}$.
        Note that this means that there is a homomorphism from~$LZ(i,\binom{t}{k-1}!)$ into~$H$ that maps the starting set to~$\{f(v^e_{11}),\dots,f(v^e_{1k-1})\}$.
        Furthermore, observe that for any~$j,j'\in\mathds{N}$ with~$j\leq j'$ there is a homomorphism from~$LZ(j',\binom{t}{k-1}!)$ to~$LZ(j,\binom{t}{k-1}!)$ that preserves the starting set.
        Therefore, there is a homomorphism from~$LZ(\binom{t}{k-1}+1,\binom{t}{k-1}!)$ to~$H$ that maps the starting set to~$\{f(v^e_{11}),\dots,f(v^e_{1k-1})\}$.
        Since the above holds for all~$e\in X^{(2)}$, we obtain a homomorphism from~$GLZ(s,\binom{t}{k-1}+1,\binom{t}{k-1}!)$ to~$H$.
        Thus,~$\ex_{\text{hom}}(t,GLZ(s,\binom{t}{k-1}+1,\binom{t}{k-1}!))\leq(\alpha_s+\varepsilon)\binom{t}{k}$.
        Since the sequence $$\frac{\ex_{\text{hom}}\Big(m,GLZ(s,\binom{t}{k-1}+1,\binom{t}{k-1}!)\Big)}{\binom{m}{k}}$$ is non-increasing in~$m$, this implies that~$\pi_{\text{hom}}(GLZ(s,\binom{t}{k-1}+1,\binom{t}{k-1}!))\leq\alpha_s+\varepsilon$.
        Thus, we have
        \begin{align*}
            \alpha_s\leq &\,\pi\Big(GLZ\big(s,\binom{t}{k-1}+1,\binom{t}{k-1}!\big)\Big)\nonumber\\
            =&\,\pi_{\text{hom}}\Big(GLZ\big(s,\binom{t}{k-1}+1,\binom{t}{k-1}!\big)\Big)\leq\alpha_s+\varepsilon\,,
        \end{align*}    
        where the first inequality holds since if~$H$ contains~$GLZ(s,\binom{t}{k-1}+1,\binom{t}{k-1}!)$, then there exists a homomorphism from~$GL(s,\ell)$ into~$H$ for all~$\ell\in\mathds{N}$.
        In fact, for all integers~$m\geq2$ and~$i\geq1$, there is a homomorphism from~$GL(s,\ell)$ into~$GLZ(s,i,m)$ for all~$\ell\in\mathds{N}$, so that~$\alpha_s\leq\pi(\{GLZ(s,i,m):\:i,m\in\mathds{N},m\geq 2\})$.
        Since in the above argument~$\varepsilon>0$ was arbitrary,~$\pi(\{GLZ(s,i,m):\:i,m\in\mathds{N},\,m\geq2\})=\alpha_s$.

    \section{Concluding remarks}

    Our earlier paper~\cite{CS:25} showed that~$\Pi^{(k)}$ contains a subset of order type~$\omega 2$ when~$k \ge 3$.
    This 
    is already enough to distinguish it from~$\Pi^{(2)}$, which has order type~$\omega$.
    In this paper, we went further, showing that~$\Pi^{(k)}$ contains a subset of order type~$\omega^2$ when~$k \ge 3$.
    This is still likely far from the truth and we conjecture that~$\Pi^{(k)}$ contains subsets of any countable order type when~$k \ge 3$.
    However, it would already be interesting to push our techniques to handle, say,~$\omega^3$,~$\omega^\omega$, or~$\epsilon_0$.
    It may also be that~$-\Pi^{(k)}$ contains subsets of any countable order type when~$k \ge 3$.
    This might be difficult, as finding a subset of order type~$\omega$ would already show that~$\Pi^{(k)}$ is not well-ordered, itself an interesting open problem.
    Finally, we note that a result of Pikhurko~\cite{P:14} saying that~$\Pi^{(k)}_{\infty}$ has the cardinality of the continuum for~$k \ge 3$ implies that~$\Pi^{(k)}_{\infty}$ and, therefore,~$\Pi^{(k)}_{\text{fin}}$ has uncountably many accumulation points.
    The same may well be true of~$\Pi^{(k)}$.

    \section*{Acknowledgements}
        We thank Forte Shinko and Mathias Schacht for interesting discussions.
        The first author was supported by NSF Awards DMS-2054452 and DMS-2348859 and by a Simons Visiting Professorship, while  
        the second author was partially supported by the Young Scientist Fellowship IBS-R029-Y7.
        This research was also supported by the National Science
        Foundation under Grant No. DMS-1928930 while the authors were in
        residence at the Simons Laufer Mathematical Sciences Institute in
        Berkeley, California during the Spring 2025 semester on Extremal Combinatorics.

\end{document}